\documentclass[11pt]{amsart}

\usepackage{amsmath,amssymb,amsfonts}
\usepackage{graphicx}
\usepackage{hyperref}

\title[CONSTRUCTION OF A CLOSED HYPERBOLIC SURFACE]{CONSTRUCTION OF A CLOSED HYPERBOLIC SURFACE OF
ARBITRARILY SMALL EIGENVALUE OF PRESCRIBED SERIAL
NUMBER}
\author{SUSOVAN PAL}

\date{April 01 ,2012.}

\begin{document}

\begin{abstract}
In this paper we construct, for given any small positive number
$\epsilon$ and given natural number $n$, and given any closed hyperbolic surface $M$, a
closed hyperbolic covering surface $\widetilde{M}$, such that its $n$-th eigenvalue is less than
$\epsilon$. An application of this result will also be discussed. The main result follows
from the techniques used in B.Randol's paper in 1974 [Ran]. Here I give a new
and geometric proof of the main result.
\end{abstract}

\maketitle

\noindent
1991 Mathematics Subject Classification. 35P15(49G05 53C20)(primary), and 58J50(secondary).

\noindent
Key words and phrases. Riemann Surfaces,Hyperbolic Geometry, Differential geometry.

\section{Introduction and Preliminaries}

A closed hyperbolic surface is a compact surface without boundary whose Gauss-
ian ( or sectional ) curvature is $-1$.

One can show that, by considering the Euler characteristic of the surface and
using Gauss-Bonnet theorem for closed surface that closed orinted surface can sup-
port a complete Riemannian metric with constant Gaussian curvature $-1$ (called
the hyperbolic metric on the surface) if and only if its genus is greater than or equal
to $2$. In this paper, we will be primarily concerned with hyperbolic metrics.

Laplacian on $M$ is the linear operator acting on the space of all smooth functions
on $M$, defined by $\Delta f = \mathrm{div}(\nabla f)$ and let $\lambda_n$ be the $n$-th eigenvalue of the Laplace
operator. For a closed oriented surface, it is known that, the spectrum, i.e. the set
of all the eigenvalus of the Laplacian is always discrete, countable and the eigen-
values are non-negative[Ch], so we can talk about the $n$-th eigenvalue of Laplacian
operator. Some authors define the Laplacian operator as $\Delta f = -\mathrm{div}(\nabla f)$ , so for
them the eigenvalues would be non-positive. Eigenvalues of the Laplace operator
has been an area of continuous study and research, because of its obvious connec-
tion with Physics and other areas in Mathematics. There have been much reserach
on the upper bound on the eigenvalues and also on whether we can produce closed
surfaces with small eigenvalues. Good references for the geometry and speactra for
closed hyperbolic surfaces are the books by P. Buser[Bu] and I. Chavel[Ch], among
others. Another good reference for both compact and non-compact Riemann sur-
faces is the book by Nicolas Bergeron[Bergeron]. I will generalize a result mentioned
in [Bergeron].

The main result, theorem 3.4, follows from the techniques involving Selberg's
trace formula used by Prof. Burton Randol's 1974 paper "Small Eigenvalues of the
Laplace Operator on Compact Riemann Surfaces"[Ran]. He also mentioned the
result in Chapter 11 of Issac Chavel's book "Eigenvalue in Riemannian Geometry".
Here I am giving an alternative proof of his result using more geometric techniques
and elementary methods.

\section{Organization of the paper}

In section 3, I will state the the minimax principle (theorem 3.1) used in the
context of eigenvalues, and the result mentioned in [Bergeron],the lemma used there
to prove that result, which is crucial in the proof of the lemma. Then in the same
section, I will state the main theorem 3.4 ,which is a generalization of the result in
[Bergeron], and the lemma 3.5 used to prove the main theorem, which generalizes
the previous lemma 3.3, and hinges on the minimax principle. I will give a complete
proof of the lemma and the main theorem.

In section 4, I will mention an application of the main theorem to a result of
R.Schoen, S. Wolpert and S. T. Yau [SWY], whose another proof is given in [DPRS].
I will state the relevant defnitions.

I also take the pleasure to thanks Prof. Feng Luo, my doctoral thesis advisor for
suggesting this problem, for various helpful discussions and encouragement in this
subject, and my colleague Tian Yang, Dr. Marius Beceanau and Ali Maalaoui for
giving me suggestions on typesetting this article.

\section{Statement and Proof of the main thorem}

We will start with the minimax principle.
Let $W^{1,2}(M)$ denote the Sobolev
space of functions on $M$ whose first order distributional derivatives exist on $M$ and
are (globally) square-integrable on $M$.

\medskip
\noindent
Theorem 3.1 ( Minimax principle ) :

\noindent
Let $f_1, f_2, ... f_{k+1}$ be continuous functions on $M$ such that they lie in the Sobolev
space $W^{1,2}(M)$ and assume that volume of ( support of $f_i \cap$support of $f_j$) $= 0$ $\forall 1 \le
j \le k+1$. Then the $k$-th eigenvalue $\lambda_k$ of $M$ satisfes the upper bound:
\[
\lambda_k \le \max_{1\le i\le k+1}\frac{\bigl(\|\nabla f_i\|_2\bigr)^2}{\bigl(\|f_i\|_2\bigr)^2}.
\]

For a proof of Minimax principle, please see Peter Buser's book [Bu] or in Issac
Chavel's book [ Ch ]. Let us remark here that some authors also use the symbol
$H^1(M)$ or $W^1(M)$ in stead of $W^{1,2}(M)$ to denote the corresponding Sobolev spaces.

Next, let us state the theorem from [Bergeron] that we intend to generalize.

\medskip
\noindent
Theorem 3.2 (Bergeron):

\noindent
Given any connected, closed, hyperbolic surface $M$, and given any $\epsilon > 0$, there
exists a finite cover $\widetilde{M}$ of $M$ such that its 1-st eigenvalue $\lambda_1(M) < \epsilon$

\medskip
To prove his theorem, [Bergeron] used the following (technical) lemma which we
will generalize as well :

\medskip
\noindent
Lemma 3.3 :

\noindent
Let $M$ be a closed hyperbolic surface such that $M = A \cup B$ where $A$ and $B$ are
two connected compact sets satisfying $A\cap B = \cup_{i=1}^n \gamma_i$, where $\gamma_i$'s are simple closed
geodesics in $M$.Let $l(\gamma_i)$ denote the length of $\gamma_i$ and let :
\[
h=\frac{\sum_{i=1}^n l(\gamma_i)}{\mathrm{minimum}\{\mathrm{area}(A), \mathrm{area}(B)\}}.
\]
Let $\eta>0$ be a positive number such that $\eta$ -neighborhood of every $\gamma_i$ is embedded
in $M$. Then there exists a constant $C(\eta)$, depending only on $\eta$ such that the first
positive eigenvalue $\lambda_1(M)$ satisfies : $\lambda_1(M)\le C(\eta)(h+h^2)$.

The proof of Theorem 3.2 [Bergeron] follows from lemma 3.3.

Finally, we state the main theorem of this paper :

\medskip
\noindent
Theorem 3.4. (Main Theorem):

\noindent
Given any connected, closed, hyperbolic surface $M$, given any natural number
$n$ ; and given any $\epsilon > 0$; there exists a finite cover
$\widetilde{M}$ of $M$ such that its $n$-th
eigenvalue $\lambda_n(\widetilde{M}) < \epsilon$ .

\medskip
Proof of the main theorem will follow from the following lemma :

\medskip
\noindent
Lemma 3.5 :

\noindent
Let $\widetilde{M}$ be a closed hyperbolic surface such that $\widetilde{M} = \cup_{i=1}^{n+1}A_i$ and
$A_1 \cap A_2 =
\gamma_1$, $A_2 \cap A_3 = \gamma_2$, ..., $A_n \cap A_{n+1} = \gamma_n$, $A_{n+1} \cap A_1 = \gamma_{n+1}$, where $A_i$'s are closed
subsets of
$\widetilde{M}$ and $\gamma_i$ 's are pairwise disjoint simple closed geodesic in
$\widetilde{M}$, and
$A_i \cap A_j = \emptyset$ $\forall j \ge i + 2$ except that $A_1 \cap A_{n+1} = \gamma_{n+1}$. Further assume that areas
of all $A_i$ and lengths of all the $\gamma_i$ 's are equal,and that $\eta$ neighborhood of each $\gamma_i$ is
embedded in $\widetilde{M}$. Then we have : $\lambda_n(\widetilde{M}) \le C(\eta)(h+h^2)$, where $C(\eta)$ is a positive
constant depending only on $\eta$, and $h = (n + 1)\cdot \dfrac{l(\gamma_1)}{\mathrm{area}(A_1)}$ .

\medskip
\noindent
Proof:

\noindent
We will use the minimax principle [Ch] to prove the lemma. We will produce
$(n + 1)$ functions $g_1, g_2, ..... g_{n+1}$ on $\widetilde{M}$ such that
\[
\frac{\bigl(\|\nabla g_i\|_2\bigr)^2}{\bigl(\|g_i\|_2\bigr)^2}\le C(\eta)(h+h^2).
\]

\noindent
Define for $t$ small positive,
\[
A_i(t)=\{z\in A_i:\mathrm{dist}(z,\gamma_i)\le t\}
\]
So $A_i(t)$ is a half-collar around the simple closed geodesic $\gamma_i$.

\noindent
Next, define the functions $f_i:\widetilde{M}\to \mathbb{R}$ by :
\[
f_i(z)=
\begin{cases}
\frac{1}{t}\,\mathrm{dist}(z,\gamma_i) & \text{if } z\in A_i(t),\\
1 & \text{if } z\in A_i\backslash A_i(t),\\
0 & \text{if } z\in \widetilde{M}\backslash A_i.
\end{cases}
\]

\noindent
Then , $\bigl(\|\nabla f_i\|_2\bigr)^2 = \dfrac{1}{t^2}.\mathrm{area}(A_i(t))$

\noindent
And, $\bigl(\|f_i\|_2\bigr)^2 \ge \mathrm{area}\bigl(A_i\backslash(A_i(t))\bigr)$

\noindent
It is clear that $f_i \in C^0(M)\cap W^{1,2}(M)\forall 1\le i\le n$. Now,
\[
\frac{\bigl(\|\nabla f_i\|_2\bigr)^2}{\bigl(\|f_i\|_2\bigr)^2}
\le \frac{1}{t^2}\cdot
\frac{\mathrm{area}(A_i(t))}{\mathrm{area}(A_i)-\mathrm{area}(A_i(t))}
\le \frac{1}{t^2}\cdot
\frac{l(\gamma_i)\sinh(t)}{\mathrm{area}(A_i)-l(\gamma_i)\sinh(t)}
\]
\[
\le \frac{1}{t^2}\cdot \frac{1}{n+1}\cdot
\frac{h.\mathrm{area}(A_i)\sinh(t)}{\mathrm{area}(A_i)-\frac{1}{n+1}.h.\mathrm{area}(A_i)\sinh(t)}
\le \frac{1}{t^2}\cdot \frac{1}{n+1}\cdot
\frac{h.\sinh(t)}{1-\frac{1}{n+1}.h.\sinh(t)}
\]
\[
\le \frac{h.\sinh(t)}{t^2(1-\sinh(t))}
\le \frac{h}{t(1-\sinh(t))}
\le \frac{2h}{\eta}\cdot \frac{1}{(1-\sinh(t))}
\le \frac{2}{\eta}\,h(1+h)
\]
\[
= C(\eta).h(1+h), \ \ \text{where } C(\eta)=\frac{2}{\eta}.
\]

\noindent
This completes the proof of the lemma 3.5.

\medskip
\noindent
Proof of the main theorem 3.4 :

\noindent
Proof of the main theorem :

\noindent
As genus of $M$ is $\ge 2$, there exists a simple closed geodesic $\gamma$ embedded in $M$ such
that $M\backslash\gamma$ is connected. Take $\eta$ positive such that $\eta$ -neighborhood of $\gamma$ is embedded
in $M$. Fix $\eta$ once and for all. Now, for each natural number $N$, construct a cover $\widetilde{M}$
of $M$ of degree $(n+1)N$ in the following way : Take $(n+1)N$ copies of $M\backslash\gamma$ and join
them in a cyclical way, i.e. each copy of $M\backslash\gamma$ is joined to two other and different
copies of $M\backslash\gamma$. Then there exists $(n+1)$ lifts of $\gamma$ cutting $\widetilde{M}$ into $(n+1)$ pieces
$A_1, A_2, ... A_{n+1}$, each one formed by $N$ fundamental domains for the action for the
covering $\widetilde{M}$ such that $A_1\cap A_2=\gamma_1$, $A_2\cap A_3=\gamma_2$, ..., $A_n\cap A_{n+1}=\gamma_n$,
$A_{n+1}\cap A_1=\gamma_{n+1}$,such that each $A_i$ is a union of $N$ copies of $M\backslash\gamma$ ,which is a disjoint union
except for a set of measure zero. Then for each $i$, $\mathrm{area}(A_i)=N.\mathrm{area}(M)$. ( see the
corresponding figure of the $(2+1).2 = 6$-fold covering surface of $M$ for $n=2, N=2$
below.)

\begin{figure}[ht]
\centering
\includegraphics[width=0.9\textwidth]{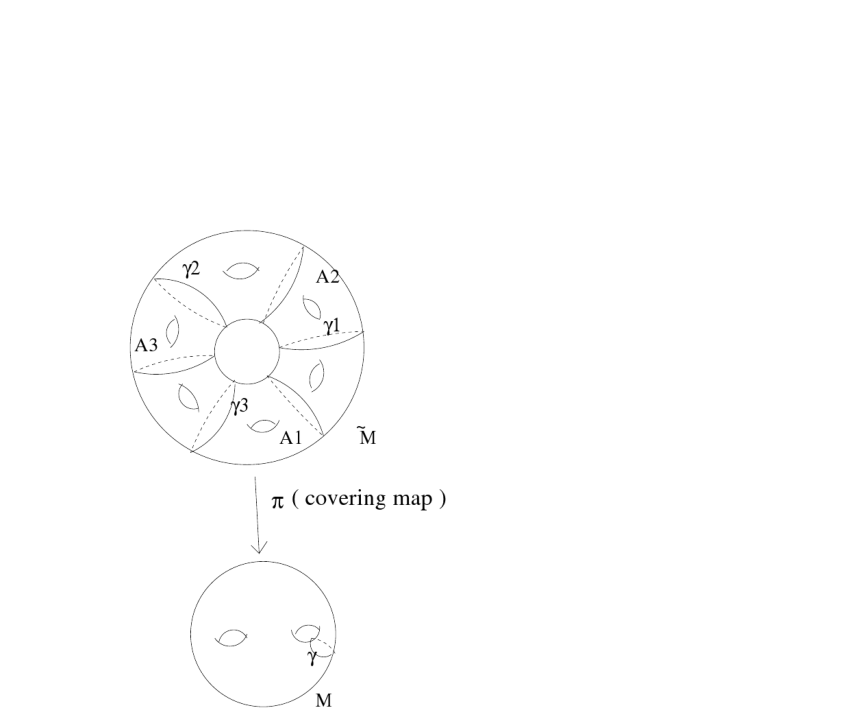}
\par\medskip
A picture of the covering surface when $n = 2, N = 2$.
\end{figure}

\noindent
Then, by the previous lemma,
\[
\lambda_n(\widetilde{M}) \le C(\eta)\cdot\left[\frac{l(\gamma)}{N.\mathrm{area}(M)}+\left(\frac{l(\gamma)}{N.\mathrm{area}(M)}\right)^2\right]\to 0
\ \ \text{as } n\to \infty
\]
This proves the main theorem of the article.
\hfill$\square$

\section{Application of the main theorem}

Let us first define some quanitities already defined in [SWY],[DPRS].

We define the quantity $l_n$ already defined in [SWY].

\medskip
\noindent
Definition 4.1 :

\noindent
Fix a closed hyperbolic surface $M$. For fixed $n$ let $C$ stand for any ( finite )
collection of simple closed geodesics in $M$ such that the complement of the union
of the geodesics in $C$ is a disjoint union of $(n+1)$ components. Let $C_n$ denote
the family of all such $C$'s.Let $l(C)$ denote the sum of lengths of all the geodesics in
$C$. Let $l_n$ denote the infimum of all $l(C)$ where $C$ varies in $C_n$ .It is easy to show
directly using the definition that, $l_n \le l_{n+1}$.

In [SWY], [DPRS], it is shown that $\lambda_n$ is related to a geometric quantity $l_n$,
where $l_n$ is the quantity defined in the definition 4.1. Then the results of the pa-
pers above show us that :

\medskip
\noindent
Theorem 4.2 ( SWY, DPRS ) :

\noindent
With $l_n$ defined as above, we have $C_1(g).l_n \le \lambda_n \le C_2(g).l_n \forall 1\le n\le (2g-3)$ ,
where $C_1(g), C_2(g)$ are constants depending on only the genus $g$ of the surface.

For a proof, see [SWY] or [DPRS]. Here we can easily prove as a
Corollary of the main theorem :

\medskip
\noindent
In the theorem 3.1 above, we cannot make $C_1(g)$ independent of $g$.

\medskip
\noindent
Proof of corollary:
For large $N$, the $n$-th eigenvalue of $\widetilde{M}$ that we just constructed
is arbitrarily close to zero but $l_n(\widetilde{M}) \ge l_1(\widetilde{M}) \ge l_1(M)$; since image of any family of
geodesics that cut $\widetilde{M}$ into two pieces cut $M$ into two pieces as well, and the image
of any geodesic in the family cutting $\widetilde{M}$ has the same length of its image, and two
geodesics in $\widetilde{M}$ could be identified in $M$ . But $l_1(M)$ is a fixed positive number
since $M$ is fixed once and for all. So $C_1(g)$ cannot be made independent of $g$ :
note that the genus of the covering surfaces go to infinity as $\epsilon$ is made arbitrarily
smaller.
\hfill$\square$

\medskip
A way to prove the dependence of $C_2(g)$ on the genus $g$ could be to construct a
sequence of hyperbolic surfaces from $M$ with $n$-th eigenvalues going to $\infty$ and their
$l_n$ being less than or equal to that of $M$.

\section*{References}

\noindent
[1] [Bergeron]Nicholas Bergeron: Le spectre des surfaces hyperboliques

\noindent
[2] [SWY] R. Schoen, S. Wolpert and S. T. Yau:Geometric bounds on the low eigenvalues of
a compact surface, Proc. Sympos. Pure Math., Vol. 36, Amer. Math. Soc, Providence, RI,
1980, pp. 279-286.

\noindent
[3] [DPRS]J. Dodziuk, T. Pignataro, B. Randol and D. Sullivan: Estimating small eigenvalues
of Riemann surfaces, Contemporary Mathematics 64 (1987), 93-121.

\noindent
[4] [Bu]Peter Buser : Geometry and Spectra of Compact Riemann Surfaces , Birkhauser, pp.
213-215.

\noindent
[5] [Ch]Issac Chavel : Eigenvalues in Riemannian Geometry : Academic Press, Chapter 1, 11.

\noindent
[6] [Ran] Burton Randol :
Small Eigenvalues of the Laplace Operator on Compact Riemann
Surfaces , Bulletin of the American Mathematical Society 80(1974)996-1000 .

\medskip
\noindent
Department of Mathematics, Rutgers University, New Brunswick, NJ 08854

\noindent
E-mail address: susovan97@gmail.com

\end{document}